\documentclass{article}
\usepackage{graphicx} 
\usepackage{subcaption}
\usepackage[style=numeric, maxbibnames=12, sorting=ynt]{biblatex}
\usepackage{tikz}
\usepackage{pgfplots}
\pgfplotsset{compat=1.18}  

\usepackage{float}

\usepackage{amsmath,amssymb,amsthm}
\usepackage[margin=1in]{geometry}
\usepackage[T1]{fontenc}
\usepackage{lmodern}

\newtheorem{proposition}{Proposition}
\newtheorem{corollary}{Corollary}
\newtheorem{remark}{Remark}
\newtheorem{theorem}{Theorem}

\newcommand{\E}{\mathbf{E}}
\newcommand{\Var}{\mathbf{Var}}
\newcommand{\Cov}{\mathbf{Cov}}
\newcommand{\Prob}{\mathbb{P}}

\title{Fluctuations in uniform one-dimensional tropical sandpile}
\author{Mikhail Shkolnikov and Alexander Varypaev}

\date{}

\begin{document}

\maketitle
\begin{abstract}
The primary study object of this note is fluctuations in a uniformly perturbed tropical sandpile on a segment. As the number of perturbations tends to infinity, the centered odometer, normalized by the square root of their number, has as its limit a Gaussian field with an explicitly computed covariance. By an elementary analytic derivation supported numerically, this covariance is a two-variable piecewise quartic function on the square, cubic in each of the variables, failing to be smooth along the diagonal. Seen as a kernel of an integral operator, it is interpreted as the Green function for a fourth-order linear differential equation describing the vertical deflection of an Euler-Bernoulli prismatic slender beam under a given transverse load. 

\end{abstract}
\section{Introduction}
The tropical sandpile model was introduced in \cite{kalinin2018self} as the conceptualization of the scaling limit announced in \cite{kalinin2016tropicalsandpiles} and established in \cite{kalinin2026tropicalmodels}. This scaling limit describes the small-mesh behavior of perturbations of the maximal stable state of the Abelian sandpile, a cellular automaton demonstrating self-organized criticality, on a convex portion of the lattice. For instance, one observable is the number of local operations per site in a relaxation, called toppling function or odometer, for small mesh it approximates the minimal non-negative tropical series (an infimum of a collection of affine-linear functions with integer gradients) with the corner locus containing the perturbation points. For the convex domain $\Omega$ and the finite set of perturbation points $P$ in the interior of $\Omega,$ we denote this tropical series $F_{\Omega,P}:\Omega\rightarrow[0,+\infty).$

Quite recently, it was shown in \cite{kalinin2026many} that in dimension two, for an infinite sequence of points $p_1,p_2\dots$ in the interior of $\Omega,$ such that $n^{-1}\sum_{k=1}^n\delta_{p_k}$ converges to a given measure $\rho$ with the support contained in the interior of $\Omega,$ the rescaled sequence of tropical odometers $n^{-1/2}F_{\Omega,\{p_1,\dots,p_n\}}$ converges to the solution of the Monge-Amp\`ere equation with the source given by the density of $\rho$ and Dirichlet boundary conditions. Earlier numerical simulations for the uniform measure case on the square \cite{kalinin2023statistics}, as well as exact solvability of the one-dimensional \cite{shkolnikov2023relaxation} case, played a significant role in the conception of this result; see \cite{lupercio2026symmetry} for some historical details and more physically-styled narration.

In this note, we initiate the study of asymptotics of fluctuations in tropical relaxation, giving their description in the simplest scenario of the uniformly-driven perturbation on an interval. Concretely, we assume that each perturbation point $p_k$ is uniformly distributed on $[0,1],$ and the points are independent. The task we are accomplishing in the next section is the description of how the random piecewise linear function $n^{-1}F_{[0,1],{p_1,\dots,p_n}}$ deviates from its deterministic limit $t\mapsto \frac{1}{2} t(1-t)$ for large $n.$

One result is that for each given point in $[0,1]$ the fluctuation is a Gaussian random variable with deviation of order $\sqrt{n}.$ Deriving this result doesn't use much apart the central limit theorem and basic structure of one-dimensional tropical sandpile. In addition, and still within standard elementary probabilistic framework we derive an explicit formula for the covariance of the fluctuations at two points of $[0,1]$, which turns out to be proportional, as soon as the number of perturbations $n$ is at least $2$, to $n+1$ times a symmetric piecewise-quartic function $C\colon[0,1]^2\rightarrow\mathbb{R}$. 

Somewhat surprisingly, apart from an explicit analytic expression, this function has a meaningful interpretation within Euler-Bernoulli beam theory. Namely, the operator sending a measure $\nu$ on $[0,1]$ to a function $f_\nu(x)=\int_0^1 C(x,y)d\nu(y)$ corresponds to a vertical displacement of a slender homogeneous beam under the load $\nu,$ where the heights of the beam at the ends of the interval are taken to be $0,$ and the first and the second derivatives at one end are equal to those at the other -- this boundary conditions allow in practical application for stacking the beam structure periodically, see Figure 3 for an example. We expect that non-uniform distributions of tropical sandpile perturbations would lead to inhomogeneous beams. 

\section{Analytic derivations}
\noindent Let $p_{1,\ldots,n} = (p_i)_{i=1}^n$ be perturbation points. From theory of one-dimensional tropical sandpile model on the interval $[0,1]$ it is known that there is a point $q$ defined by finite relaxation procedure which can be calculated as follows (see \cite{shkolnikov2023relaxation})
\[
1 + \left\lfloor \sum_{i=1}^{n} p_i \right\rfloor - \sum_{i=1}^{n} p_i.
\]
Thus, $q$ is the distance from the sum of the coordinates to the nearest integer above it. This is not some new, randomly placed point, rather, it is determined by the original points. Thus the finite slopes of $g_n$ are changed at $p_1,\ldots,p_n, q$, and in each one, the slope decreases by one. Let $v_1 \le v_2 \le \ldots \le v_{n+1},$  be the ordered locations of slope changes of the odometer $g_n.$ 

\begin{proposition}
The tropical function $g_n$ can be written as
\begin{equation}
g_n(t) = \left( n+1 - \sum_{i=1}^{n+1} v_i \right) t - \sum_{i=1}^{n+1} (t - v_i)_+, \quad \text{where} \quad (t - v_i)_+ = \begin{cases} 0, & t \le v_i, \\ t - v_i, & t > v_i .\end{cases} \tag{1}
\end{equation}
\end{proposition}

\begin{proof}
Denote by $a$ the initial slope of the odometer, then due to fact that function starts at $(0,0)$ then before $v_1$ i.\,e.\ for $0 \le t \le v_1$ we have $g_n(t) = at$. For $v_1 \le t \le v_2$ we have slope $a-1$ and in order to write this fact taking into account the continuity of the function $g_n(t)$ we should subtract $t - v_1$ but only for $t \ge v_1$, i.\,e.\ for $0 \le t \le v_2$ we will have $g_n(t) = at - (t - v_1)_+$. Each subsequent slope is accounted for in exactly the same way, thus we get:
\begin{equation}
g_n(t) = at - \sum_{i=1}^{n+1} (t - v_i)_+ \tag{2}.
\end{equation}
Here we already ensured that $g_n(0) = 0$. In order to define initial slope $a$ we should also verify that $g_n(1) = 0$ and from this condition we have
\begin{equation}
0 = a - \sum_{i=1}^{n+1} (1 - v_i) \;\Rightarrow\; a = n+1 - \sum_{i=1}^{n+1} v_i. \tag{3}
\end{equation}
\end{proof}

Substituting (3) in (2) we finally get (1).

Expression (1) is very useful for further exploration of probability properties of $g_{n+1}(t)$.

\begin{proposition}
Let $(p_i)_{i=1}^n$ be a random set of independent random points and each $p_i$ -- uniformly distributed, i.\,e.\ $p_i \sim U[0,1]$ then the following convergence holds
\begin{equation}
\frac{g_n(t) - \frac{n+1}{2}t(1-t)}{\sqrt{n}} \xrightarrow[n \to \infty]{} N \left( 0, \frac{t^2(1-t)^2}{12} \right).\tag{4}
\end{equation}
\end{proposition}

\begin{proof}
We can rewrite \(g_n(t)\) as follows:
\begin{equation}
g_n(t)=\sum_{i=1}^{n}X_i(t)+R_n(t),\tag{5}
\end{equation}
where
\[
X_i(t)=t(1-p_i)-(t-p_i)_+=
\begin{cases}
t(1-p_i), & t\le p_i,\\
p_i(1-t), & t>p_i,
\end{cases}
\qquad
R_n(t)=t(1-q)-(t-q)_+.
\]
Now we are ready to calculate mathematical expectation of \(g_n(t)\):
\begin{equation}
\mathbf{E}X_i(t)
=\int_0^t (1-t)x_i\,dx_i+\int_t^1 t(1-x_i)\,dx_i
=\frac{t(1-t)}{2}=\mu(t).\tag{6}
\end{equation}
It can be shown that \(q\sim U[0,1]\), then
\[
\mathbf{E}g_n(t)
=\sum_{i=1}^{n}\mathbf{E}X_i(t)+\mathbf{E}R_n(t)
=n\mu(t)+\mu(t)
=\frac{(n+1)}{2}t(1-t).
\]
Let's calculate second moment for one term of sum (5)
\[
\mathbf{E}X_i^2(t)
=\int_0^t (1-t)^2x_i^2\,dx_i+\int_t^1 t^2(1-x_i)^2\,dx_i
=\frac{t^2(1-t)^2}{3}.
\]
The variance of \(X_i(t)\) is determined by subtracting the square of mathematical expectation from second moment:
\[
\operatorname{Var}X_i(t)
=\mathbf{E}X_i^2(t)-(\mathbf{E}X_i(t))^2
=\frac{t^2(1-t)^2}{3}-\frac{t^2(1-t)^2}{4}
=\frac{t^2(1-t)^2}{12}
=(\sigma(t))^2.
\]
Let's center function \(g_n(t)\) under \(\sqrt{n}\) normalization, i.e.
\[
Y_n(t)
=\frac{g_n(t)-\mathbf{E}g_n(t)}{\sqrt{n}}
=\frac{\sum_{i=1}^{n}(X_i(t)-\mu(t))+(R_n(t)-\mu(t))}{\sqrt{n}},
\]
\[
\operatorname{Var}Y_n(t)
=\left(1+\frac{1}{n}\right)\frac{t^2(1-t)^2}{12}
\xrightarrow{n\to\infty}\frac{t^2(1-t)^2}{12},
\qquad
\mathbf{E}Y_n(t)=0,\quad n\ge 2.
\]

The classical i.i.d.\ central limit theorem applies only to set of variables \(X_1(t),\ldots,X_n(t)\). But if we include \(R_n(t)\) into this set then we get only pairwise independence (weaker condition), because \(q_n\) is determined through \((p_i)_{i=1}^{n}\). Term \(Y_n(t)\) can be rewritten as
\[
Y_n(t)
=\underbrace{\frac{1}{\sqrt{n}}\sum_{i=1}^{n}\bigl(X_i(t)-\mu(t)\bigr)}_{A_n(t)}
+\underbrace{\frac{1}{\sqrt{n}}\bigl(R_n(t)-\mu(t)\bigr)}_{B_n(t)}.
\]
The central limit theorem can be applied to the term \(A_n(t)\):
\[
A_n(t)\xrightarrow{d}N\left(0,\frac{t^2(1-t)^2}{12}\right).\tag{7}
\]
Since \(B_n(t)\xrightarrow{p}0\), Slutsky's lemma establishes the asymptotic normality of \(Y_n(t)\).
\end{proof}
\begin{proposition}
For any $t_1,t_2 \in [0,1]$ and $n\geq 2$, covariance of two random variables $g_n(t_1)$ and $g_n(t_2)$ can be represented as follows
\[
\begin{aligned}
\Cov(g_n(t_1), g_n(t_2)) &= \frac{(n+1)s(1-u)}{12} \left( u - 2u^2 - 2s^2 + 3su \right), \\
s &= \min(t_1, t_2), \quad u = \max(t_1, t_2).
\end{aligned}
\]
\end{proposition}

\begin{proof}
Denote
\begin{equation}
C(t_1, t_2) = \Cov(X_i(t_1), X_i(t_2)) = \E X_i(t_1)X_i(t_2) - \E X_i(t_1)\E X_i(t_2). \tag{8}
\end{equation}
We start from calculation of function $C(t_1, t_2)$. From expression (6) we already know
\[
\E X_i(s) = \frac{s(1-s)}{2}; \quad \E X_i(u) = \frac{u(1-u)}{2}; \quad s \le u.
\]
It is left to calculate mathematical expectation of multiplication $X_i(t_1)X_i(t_2)$, and here we should take into account three cases:
\begin{enumerate}
\item $x_i$ is on the left side from both $t_1$ and $t_2$: $x_i \le s$
\item $x_i$ is between $t_1$ and $t_2$: $s \le x_i \le u$
\item $x_i$ is on the right side from both $t_1$ and $t_2$: $x_i \ge u$
\end{enumerate}
And thus we have
\begin{equation}
\begin{aligned}
\E X_i(t_1)X_i(t_2) &= \int_{0}^{s} (1-s)(1-u)x_i^2 \, dx_i + \int_{s}^{u} s(1-u)x_i(1-x_i) \, dx_i + {} \\
&\quad + \int_{u}^{1} su(1-x_i)^2 \, dx_i = \frac{s(1-u)}{6}\left(2u - u^2 - s^2\right).
\end{aligned} \tag{9}
\end{equation}
Subtracting from (9) the product $\E X_i(t_1)\E X_i(t_2)$ according to (8) we obtain $C(t_1, t_2)$:
\begin{equation}
C(t_1, t_2) = \E X_i(t_1)X_i(t_2) - \frac{s(1-s)u(1-u)}{4} = \frac{s(1-u)}{12}\left(u - 2u^2 - 2s^2 + 3su\right). \tag{10}
\end{equation}
We should notice the following
\begin{enumerate}
\item The contributions of different starting points $p_i$, $p_j$ have zero covariance.
\item Under $n \ge 2$ the point $q$ is independent of each individual $p_i$ that's why mixed covariance here is equal to zero.
\item Each point $p_i$, $i = 1,\dots,n$ has an impact equal to $C(t_1, t_2)$
\item Point $q$ gives the same impact as any point $p_i$, i.\,e.\ $C(t_1, t_2)$
\end{enumerate}
Therefore we can conclude that
\begin{equation}
\Cov(g_n(t_1), g_n(t_2)) = \sum_{i=1}^{n} C(t_1, t_2) + C(t_1, t_2) = (n+1)C(t_1, t_2). \tag{11}
\end{equation}
Unification of formulas (10) and (11) proves the proposition.
\end{proof}

\begin{corollary}
If $t_1 = t_2$ then $s = u$ and $\Cov(g_n(t_1), g_n(t_1)) = \Var g_n(t_1) = \frac{n+1}{12} t_1^2 (1-t_1)^2$ -- that coincides with the one in formula (4).
\end{corollary}

\begin{corollary}
There is a multidimensional central limit theorem that gives us probability distribution of vector
\[
\begin{pmatrix} Y_n(t_1) \\ Y_n(t_2) \end{pmatrix} \xrightarrow[n \to \infty]{d} N_2 \left( \begin{pmatrix} 0 \\ 0 \end{pmatrix}, \begin{pmatrix} C(t_1, t_1) & C(t_1, t_2) \\ C(t_1, t_2) & C(t_2, t_2) \end{pmatrix} \right).
\]
\end{corollary}

Now we would like to justify the Gaussian field claim for the asymptotic behaviour of odometer fluctuations. On the one hand, for the case of the interval as the underlying domain of the tropical sandpile, it would be more conventional to use the term ``Gaussian process''. However, the physical flavour of this interval is that it is not temporal but spatial; in particular, the present note serves as a conceptual preparation for the analogous study in higher dimensions, which is expected to have additional number-theoretic and geometric nuances and complications. 

On the other hand, so far we dealt with the pointwise convergence of the rescaled fluctuations. To use the word ``field'' we would need to specify the topology on the space of functions on the interval. Here, the natural choice is the one given by the sup norm. What we need to observe in addition to the above propositions is that $X_i$ and $\mu$ satisfy the 1-Lipschitz property. Together with continuity of the kernel $C,$ and by Kolmogorov–Chentsov theorem, we have the following.

\begin{theorem} For $n$-point uniform perturbation odometer $g_n$ on $[0,1],$ the sequence of random functions $$n^{-1/2}(g_n(t)-\frac{(n+1)t(1-t)}{2})$$ converges weakly in sup norm to the centered Gaussian field with covariance kernel $C$ given by (10). 
\end{theorem}

\begin{remark}
The above covariance kernel $C,$ in fact, is the Green function for the fourth-order differential equation $f''''=\sigma$ with boundary conditions, $f(0)=f(1)=0,\ f'(0)=f'(1),\ f''(0)=f''(1),$ which is the Euler–Bernoulli beam equation for a homogeneous prismatic slender beam under the load $\sigma$. See, for instance, \cite{simone2011slender} for further details, and Figure \ref{fig:beam-periodic} for graphical representation of a concrete solution.
\end{remark}

We expect that a similar result holds for non-uniform perturbation densities, resulting in a non-constant bending stiffness factor. This, and the analysis of the deformed one-dimensional tropical sandpiles of the type \cite{varypaev2025probability} will be done in the forthcoming work by the authors. The higher-dimensional case presents a promising fusion of probabilities, differential equations, and  number theory, as the fluctuation should present the non-equiaffine correction to the lattice-unimodular rescaled odometer near the limit.

\section{Numerical experiments}
From the perspective of the development of the emerging theory to which the present paper belongs, the idea that one can study fluctuations of the rescaled odometer corresponding to the stochastic perturbations, and the empirical discovery that those fluctuations form a Gaussian field, came way earlier than the analytic derivation of the previous section. Thus, these two epistemological streams are not subject to a strict axiological hierarchy, but rather the components which are allowing for the progressive increment of our understanding of tropical sandpile. 

Two Monte Carlo experiments were performed for the classical one-dimensional tropical sandpile on $[0,1]$, starting from the zero function. The first experiment examines uniform convergence of the normalized final profile to a deterministic parabola. The second compares the empirical covariance of the centered fluctuations with the analytical covariance kernel.

\subsection*{1. Uniform convergence of the normalized profile}

For each realization $r$, $n$ independent perturbation points were sampled from the uniform distribution on $(0,1)$. Let $g_n^{(r)}(t)$ denote the final tropical function for this finite set of points. The normalized profile and its deterministic limit are
\begin{equation}
h_n^{(r)}(t) = \frac{g_n^{(r)}(t)}{n}, \quad h(t) = \frac{t(1-t)}{2}. \tag{12}
\end{equation}
For each sample size, the error was measured over the entire interval:
\begin{equation}
E_n^{(r)} = \sup_{0 \le t \le 1} \left| h_n^{(r)}(t) - h(t) \right|. \tag{13}
\end{equation}
The quantity estimated in the experiment was the probability that this uniform error exceeds a fixed tolerance:
\begin{equation}
p_n(\epsilon) = \Prob\{E_n > \epsilon\}, \quad \hat{p}_n(\epsilon) = \frac{1}{M} \sum_{r=1}^{M} \mathbf{1}\{E_n^{(r)} > \epsilon\}. \tag{14}
\end{equation}
The ensemble size was $M = 1000$. Sample sizes were $n = 5 + 20k$, $k = 0, \ldots, 5000$, giving $5001$ values up to $n = 100005$. Within each realization, larger samples retained all previously generated points; different realizations used independent random samples. The three tolerances shown in Figure~1 were
\begin{equation}
\epsilon_1 = 10^{-3}, \quad \epsilon_2 = 7 \cdot 10^{-4}, \quad \epsilon_3 = 5 \cdot 10^{-4}. \tag{15}
\end{equation}

All three curves show a decreasing overall trend. Smaller tolerances yield larger exceedance probabilities, as required by the nesting of the corresponding events. The largest tolerance reaches the vicinity of zero earlier, whereas the smallest tolerance retains a visible exceedance probability near the end of the simulated range. Local fluctuations are compatible with finite-ensemble estimation and with the non-monotone evolution of individual errors.

These results illustrate uniform convergence in probability: for each fixed positive tolerance, the probability of an error exceeding that tolerance anywhere on the interval tends to zero. The experiment provides numerical evidence for this behavior over the simulated range; it does not replace the analytical convergence argument.

\begin{figure}[H]
\centering
\includegraphics[width=0.84\textwidth]{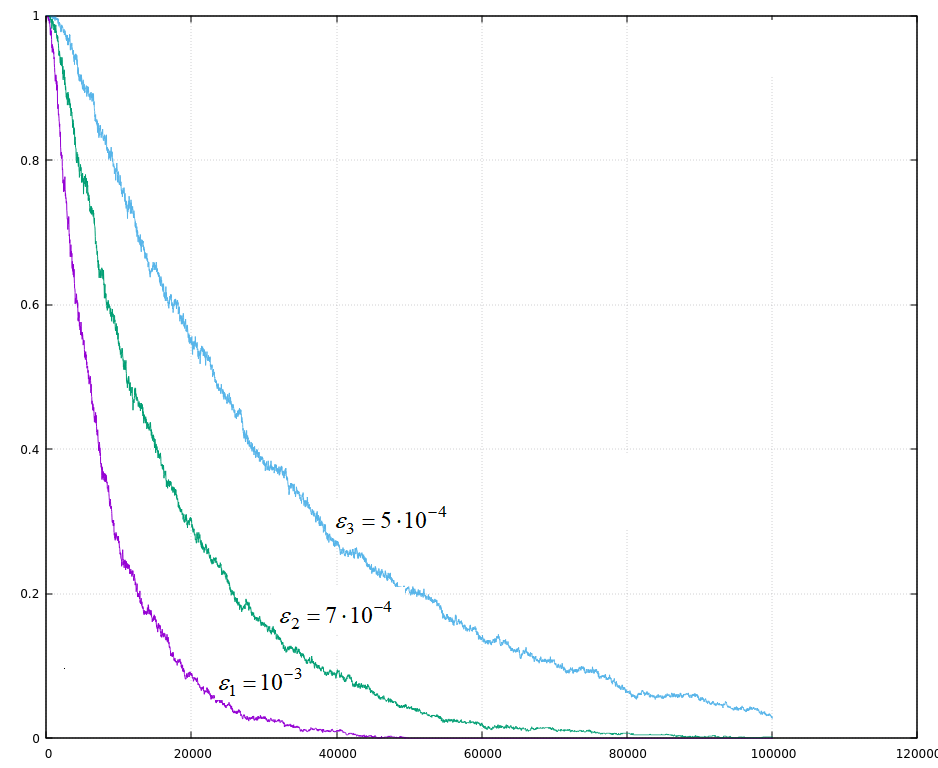}
\caption{Monte Carlo estimates of $p_n(\epsilon)$ from 1000 realizations. The horizontal axis is the number $n$ of perturbation points; the vertical axis is the estimated exceedance probability. Purple: $\epsilon = 10^{-3}$; green: $\epsilon = 7 \cdot 10^{-4}$; blue: $\epsilon = 5 \cdot 10^{-4}$.}
\label{fig:exceedance}
\end{figure}
\begin{figure}[htbp]
    \centering
    \begin{tikzpicture}
    \begin{axis}[
        width=\linewidth,
        height=6.5cm,
        xmin=-1, xmax=3,
        ymin=-0.22, ymax=0.06,
        xlabel={$x$},
        ylabel={$54\,f(x)$},
        axis lines=middle,
        samples=160,
        tick label style={font=\small},
        label style={font=\small},
        legend pos=south east,
        legend style={at={(1,0)}},
        ylabel style={xshift=-3.7em},
        xlabel style={yshift=-1.1em}
    ]

    \pgfplotsinvokeforeach{-1,0,1,2}{
        \addplot[
            blue,
            very thick,
            domain=#1:#1+0.25,
        ] {54*(((x-(#1))^3)/8
              - 3*((x-(#1))^2)/64
              - (x-(#1))/128)};

        \addplot[
            blue,
            very thick,
            domain=#1+0.25:#1+1,
        ] {54*(-((x-(#1))^3)/24
              + 5*((x-(#1))^2)/64
              - 5*(x-(#1))/128
              + 1/384)};
    }
    \addlegendimage{blue,very thick}
    \addlegendentry{periodic deflection}

    \pgfplotsinvokeforeach{-1,0,1,2}{
        \draw[->, very thick, red]
            (axis cs:{#1+0.25},0.035) --
            (axis cs:{#1+0.25},-0.15);
    }

    \addplot[only marks, mark=*, mark size=1.5pt, black]
        coordinates {(-1,0) (0,0) (1,0) (2,0) (3,0)};

    \end{axis}
    \end{tikzpicture}
    \caption{Periodic extension (four periods) of the downward-loaded
    Euler--Bernoulli beam deflection with periodic boundary conditions
    on each unit cell. Loads act at \(x=k+\tfrac14\),
    \(k=-1,0,1,2\), which is $-C(x,1/4)$ extended periodically. The vertical axis is \(54\,f(x)\) for visibility;
    the unscaled deflection is of order \(10^{-3}\). The boundary conditions ensure no kinks at the integers, as well as same curvatures to the sides, but force further bending the slender beam to the central sides of the atomic downward loads shown by red arrows. }
    \label{fig:beam-periodic}
\end{figure}
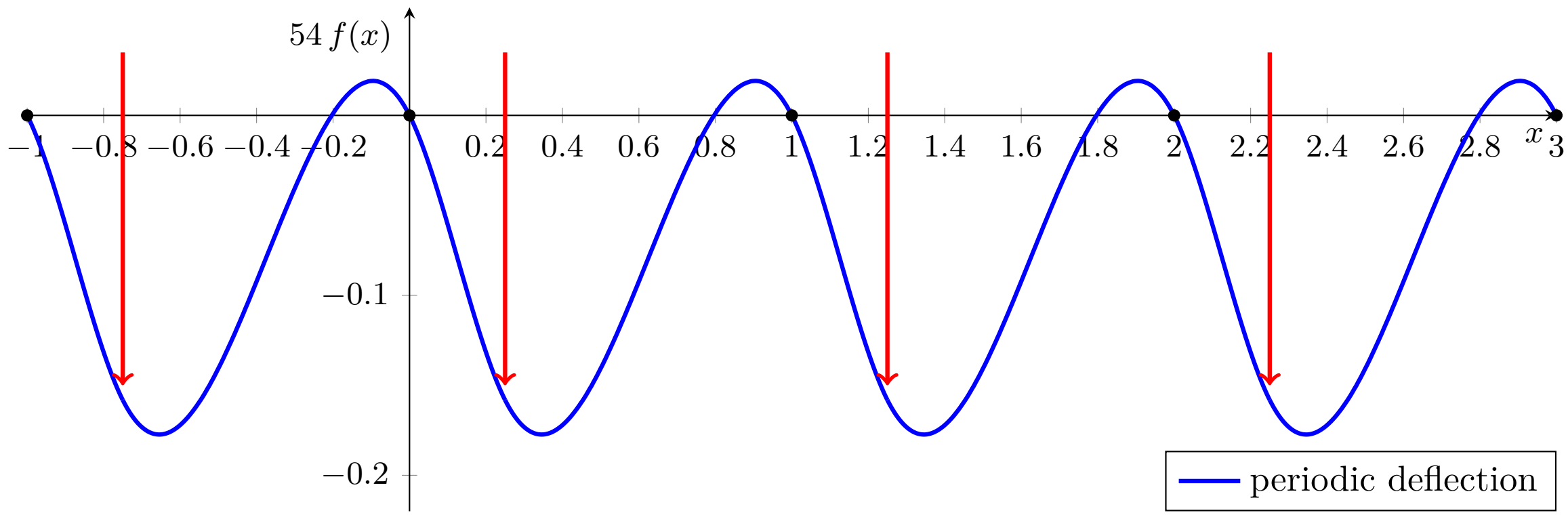

\begin{figure}[H]
\centering
\vspace{-72pt}
\begin{subfigure}{0.84\textwidth}
\centering
\includegraphics[width=\textwidth]{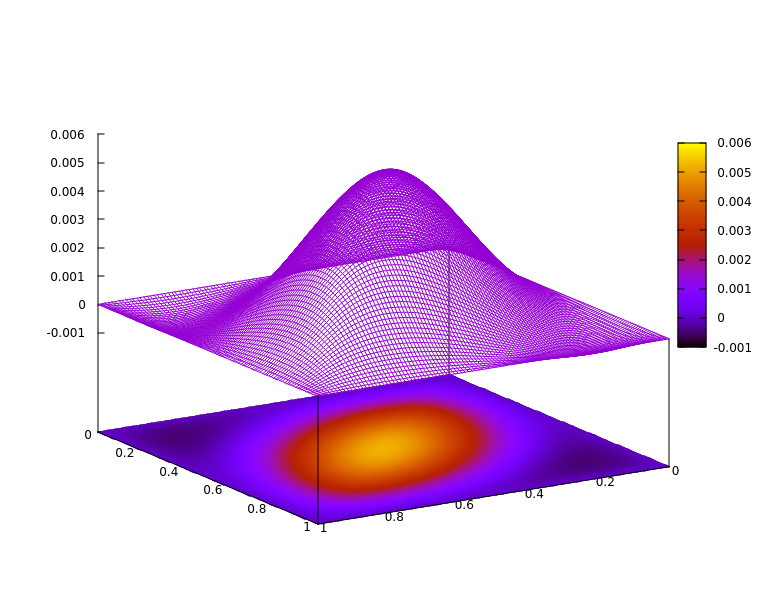}
\caption{Empirical covariance. Estimate (17), computed from 1000 realizations with $n = 10000$ and spatial step $0.02$.}
\label{fig:empirical_cov}
\end{subfigure}\hfill
\begin{subfigure}{0.84\textwidth}
\centering
\includegraphics[width=\textwidth]{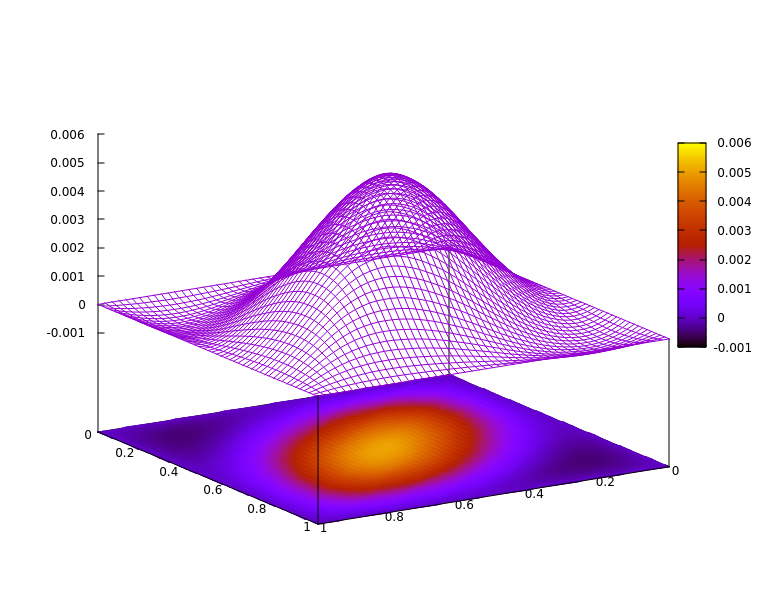}
\caption{Analytical covariance. Limiting kernel (18), evaluated on the unit square.}
\label{fig:analytical_cov}
\end{subfigure}
\caption{Three-dimensional covariance surfaces with color projections on the base plane. The two horizontal coordinates are $t_1$ and $t_2$; height and color encode covariance. The supplied figures use the same displayed color-scale limits, from $-0.001$ to $0.006$.}
\label{fig:covariance_comparison}
\end{figure}

\subsection*{2. Empirical covariance of the fluctuations}

A separate ensemble of $M = 1000$ realizations was generated with $n = 10000$ independent uniform perturbation points in each realization. The final profiles were evaluated at $t_j = 0.02j$, $j = 0, \ldots, 50$. To compare with the limiting covariance kernel, the profiles were centered by their exact expectation and normalized by the square root of $n$:
\begin{equation}
Z_n^{(r)}(t) = \frac{g_n^{(r)}(t) - \E g_n(t)}{\sqrt{n}}, \quad \E g_n(t) = \frac{n+1}{2}\,t(1-t). \tag{16}
\end{equation}
Since the exact ensemble mean of $Z_n(t)$ is zero, the covariance was estimated by the average of the outer products:
\begin{equation}
\hat{C}_n(t_j, t_k) = \frac{1}{M} \sum_{r=1}^{M} Z_n^{(r)}(t_j)\, Z_n^{(r)}(t_k). \tag{17}
\end{equation}
Each outer product uses two coordinates of the same realization, preserving their spatial dependence. Division by $M$ is appropriate because the centering uses the known expectation, rather than an estimated sample mean. The resulting $51 \times 51$ array was exported as triples $(t_j, t_k, \text{estimated covariance})$, with separate grid rows, for surface and heat-map visualization in gnuplot.

For the analytical comparison, set $s = \min(t_1, t_2)$ and $u = \max(t_1, t_2)$. The limiting covariance kernel is
\begin{equation}
C(t_1, t_2) = \frac{s(1-u)}{12}\left( u - 2u^2 - 2s^2 + 3su \right). \tag{18}
\end{equation}

The exact finite-sample covariance of the normalized fluctuations satisfies, for $n$ at least 2,
\begin{equation}
\Cov(Z_n(t_1), Z_n(t_2)) = \left( 1 + \frac{1}{n} \right) C(t_1, t_2). \tag{19}
\end{equation}
Consequently, the expectation of the estimator in (17) differs from the limiting kernel by the multiplicative factor $1.0001$ at $n = 10000$. The covariance of the unnormalized profiles is instead $(n + 1)C(t_1, t_2)$; it is not the quantity plotted here.

Figure~3 compares the supplied empirical and analytical surfaces and their projections onto the base plane. Both display a positive central peak, negative covariance regions away from the diagonal, symmetry under interchange of the two coordinates, and zero covariance at the boundary. The analytical diagonal is $C(t,t) = t^2(1-t)^2/12$, attaining its maximum $1/192$ at $t = 1/2$.

The similarity of the two plots supports the analytical covariance formula. Their remaining differences are consistent with finite-ensemble Monte Carlo variation. The comparison here is visual: no numerical discrepancy norm is inferred from the supplied images.

\section{Few more remarks}
Despite Remark 1 being a punchline of the present paper, allowing one to easily communicate the finding without reproducing a still not too messy formula for the covariance kernel of our Gaussian field, the connection with elasticity theory at the moment is shallow  -- our incompetence in that branch of human creativity doesn't allow for a deeper presentation of this barely mentioned aspect we arrived to. Moreover, the differential equation for the slender beam, and more importantly the boundary conditions of being level on the two sides and periodic in slope and curvature, is something that might not necessarily carry this particular engineering interpretation.

What we intend to do in the next expanded note is to treat the non-uniform perturbation profile with more general boundary conditions for the odometer. At this point, we expect to have a richer correspondence between tropical sandpile covariances and Green functions of Euler-Bernoulli beam equations. However, even the perfect matching at the level of formulas would not be a fully satisfactory result -- what one rather needs is a conceptual explanation for this phenomenal correspondence -- think, for instance, of the concept of Gaussian free field having the Green function of a Laplacian as its covariance kernel and all the corpus of corresponding literature across disciplines. 

At this moment we have to add that there have been a number of not yet publicly available parallel developments on the subject, which allow one to think that one can go way beyond the second order elliptic differential equations at the level of mean values, with already non-linear Monge-Amp\`ere equation starting dimension two, and fourth order equations representing covariances, whose shape is presently a mystery outside of the one-dimensional case. Tropical sandpile model itself is flexible enough for abundant numerical and hermeneutic applications, and we are only now starting to understand which of its seemingly rigid structural components represent hidden parameters subject to intentional deformations.

\section*{Acknowledgements}
Research of M.S. is funded in part by the Simons Foundation, grant [SFI-MPS-T-Institutes-00007697] and
the Ministry of Education and Science of the Republic of Bulgaria,
grant [DO1-239/10.12.2024]. The authors are grateful to Nikita Kalinin, Ernesto Lupercio and Higinio Serrano for encouragement and valuable comments, as well as to IMI-BAS and its administration for providing excellent work conditions.

\printbibliography

{\centering
\vspace{27pt}
Institute of Mathematics and Informatics \\
Bulgarian Academy of Sciences\\
Akad. G. Bonchev 8, Sofia 1113, Bulgaria \\\vspace{12pt} 

m.shkolnikov[at]math.bas.bg\\ \vspace{10pt}
varipaev[at]math.bas.bg

}

\end{document}